\documentclass[11pt]{article}
\usepackage{amsmath}
\usepackage{amsfonts}
\usepackage{latexsym}
\usepackage{verbatim}
\usepackage{epsfig}
\newtheorem{theorem}{Theorem}
\newtheorem{corollary}{Corollary}

\newtheorem{proposition}{Proposition}
\newcounter{INDEX}
\newcommand{\noi}{{\noindent}}

\newcommand{\C}{{\mathcal C}}
\newcommand{\cN}{{\mathcal N}}
\newcommand{\cR}{{\mathcal R}}
\newcommand{\cZ}{{\mathcal Z}}
\newcommand{\D}{{\mathcal D}}

\newcommand{\omb}{{\overline{\omega}}}
\newcommand{\R}{{\mathbb R}}
\newcommand{\Rn}{{{\mathbb R}^n}}
\newcommand{\Sno}{{{\mathbb S}^{n-1}}}
\newcommand{\spalpha}{{\rm span}\{\alpha\}}
\newcommand{\T}{{\mathbb T}}
\newcommand{\Tn}{{{\mathbb T}^n}}
\newcommand{\Z}{{\mathbb Z}}
\newcommand{\Zn}{{\mathbb Z}^n}
\newcommand{\Zno}{{\mathbb Z}^n\backslash\{0\}}
\begin{document}

 \title{Filling times for linear flow on the torus \\
 with truncated Diophantine conditions:\\
 a brief review and new proof}
 
% \author{H.~Scott Dumas \ and \ St\'ephane Fischler} 

\date{}

\maketitle

\vskip -0.5 true in\hskip -0.2 true in
{\begin{tabular}{ll}
 H.~Scott Dumas  &St\'ephane Fischler  \\
\small Department of Mathematics and Statistics \ \ & \small Universit\'e Paris-Saclay \\
\small University of New Mexico  & \small CNRS \\ 
\small Albuquerque, NM 87131 USA  \ \& & \small  Laboratoire de math\'ematiques d'Orsay \\
\small Department of Mathematical Sciences & \small 91405 Orsay, France \\
\small University of Cincinnati & {\tt stephane.fischler@universite-paris-saclay.fr} \\
\small Cincinnati, OH 45221 USA & \\
{\tt dumashs@ucmail.uc.edu} & \\
\end{tabular}}

\begin{abstract}  % approximately 80 words
{
We show that the geometry-of-numbers method used by A.~Bounemoura to obtain filling times for linear flow
on the torus satisfying Diophantine conditions may be extended to the case of linear flow with {\it truncated} Diophantine
conditions, and we use these methods to recover the optimal estimate first obtained by M. Berti, L. Biasco, and P. Bolle
in 2003.  We also briefly review the dynamics of linear flow on the torus, previous results, optimality, and applications of these estimates.
}
\end{abstract}


\section{Introduction}

Linear flow on the $n$-torus $\Tn$ occurs routinely in integrable and nearly integrable Hamiltonian systems, as well as in other mathematical settings.
When the direction of such flow satisfies Diophantine conditions (as happens for example on KAM tori), for given $\delta>0$ it can be shown that each orbit of the
flow becomes $\delta$-dense on the torus after a time $T$ that may be estimated in terms of $\delta$ and the Diophantine parameters.  In fact, this ``filling" occurs
even when the flow only satisfies truncated Diophantine conditions (i.e., satisfies Diophantine conditions only up to a certain critical order $N^*$).
The first crude estimates of the filling time\footnote{In earlier work, we used the term ``ergodization time" rather than ``filling time."} 
$T$ were found in the late 1980s 
by one of us (HSD) in the context of an application to physics \cite{D88}.  Since then, estimates have been refined and improved by a number of authors, and in
2003 the optimal estimate was proved for the more general case (truncated Diophantine conditions) by M.~Berti, L.~Biasco, and P.~Bolle \cite{BBB03}.  In this
paper, we give a new proof of this estimate by modifying the geometry-of-numbers method used by A.~Bounemoura  \cite{B16} to get filling times for the less
general case  (untruncated Diophantine conditions).  Bounemoura's method is in turn based on techniques developed earlier with one of us (SF) in \cite{BF13}.

In addition to our  proof, we provide background material and an overview of the filling-time problem and how it has evolved.  We hope this will make the subject 
accessible to a wider audience and draw attention to its many contributors and surprising number of applications.

The remainder of this paper is organized as follows.  In \S 2, we set out our notation and terminology, defining the filling property for linear flow on $\Tn$ and Diophantine
sets of vectors both with and without truncation.  In \S 3, we state our basic problem and main result (Theorem 1) along with corollaries for comparison with other results.
We next look at properties of Diophantine sets and discuss the dynamics of linear flow in \S 4.  Section 5 serves as a short review by surveying past results on filling times,
their optimality, and their use in applications.  In \S 6 we prove Theorem~1 by means of a proposition adapted from \cite{BF13} to treat the case of truncated Diophantine
conditions, and we briefly discuss our new proposition and overall proof.
Finally, to make the paper self contained, basic elements from geometry of numbers are presented in an appendix at the end.

\section{Notation and Terminology}

For integer dimension $n\ge2$ and $x=(x_1,\dots,x_n)\in\Rn$, $k=(k_1,\ldots,k_n)\in\Zn$, we use the Euclidean norms 
$\|x\|=(x_1^2+\cdots+x_n^2)^{1/2}$ and $\|k\|=(k_1^2+\cdots+k_n^2)^{1/2}$.
We denote by $\Tn\equiv\Rn/\Zn$ the flat $n$-torus, on which we use modular (mod $\Zn$) arithmetic (i.e., mod 1 arithmetic in each coordinate).
For $\alpha\in\Sno = \{\alpha\in\Rn\,\big| \,\|\alpha\|=1\}$, $t\in\R$, and $\theta\in\Tn$, we use the unconventional but convenient notation \ 
$\alpha_t:\Tn\to\Tn$, \ $\alpha_t(\theta) = \theta + t\alpha$ \ to denote {\sl linear flow} on $\Tn$ with unit speed and {\sl direction vector} $\alpha$.
In this paper we consider $\Zn$ to be included in $\Rn$, so when we say that a set of integer vectors is linearly independent,
or we indicate the span of a set of integer vectors, these have their ordinary meanings in $\Rn$.

Our key concept is what we call the filling property of linear flow on $\Tn$, defined as follows.

\medskip
\noi{\bf Filling $\Tn$ to within $\delta$ after time $T$.} \ 
{\sl Given $\delta\in(0,\frac{1}{2})$, we say $\alpha_t$ fills $\Tn$ to within $\delta$ after time $T$ if,  for any starting point $\theta\in\Tn$,   the
orbit segment $\{\alpha_t(\theta), 0\leq t\leq T\}$  forms a $\delta$-dense subset of $\Tn$.}
(This means that every closed ball of radius $\delta$ in $\Tn$ contains a point of the orbit segment.)

It is not difficult to see that the filling property and filling time $T$ are independent of the starting point $\theta$.  For more details about this fact, see
the introductory parts of \cite{B16} or \cite{D91}. 

\medskip
\noi{\bf Diophantine  vectors.} \ The filling property of the flow $\alpha_t$ depends strongly on the Diophantine properties of the direction vector $\alpha$.
We use the following sets.

For  dimension $n\ge2$ and parameters $\tau>n-1$, $\gamma>0$, we define $\D_n(\tau,\gamma)$, the set of {\sl Diophantine vectors} in $\Rn$ by

$\D_n(\tau,\gamma) = \bigl\{\alpha\in\Rn\,\big|\ |k\cdot\alpha|\ge\gamma\|k\|^{-\tau} \ {\rm for\ any} \ k\in\Zno\bigr\}$.

With $n, \tau, \gamma$ as above, we adjoin the parameter $N\ge1$ and define the set of {\sl truncated Diophantine vectors} by

$\D_n(\tau,\gamma,N) = \bigl\{\alpha\in\Rn\,\big|\ |k\cdot\alpha|\ge\gamma\|k\|^{-\tau} \ {\rm for\ any} \ k\in\Zn \ {\rm with} \ 0<\|k\|\le N\bigr\}$.

Finally, we attach the superscript 1 to these sets to indicate the restriction to vectors of unit length.  In other words,
$\D^1_n(\tau,\gamma) = \Sno\cap\D_n(\tau,\gamma)$ and $\D^1_n(\tau,\gamma,N) = \Sno\cap\D_n(\tau,\gamma,N)$.
We loosely refer to $\alpha$ in $\D_n(\tau,\gamma)$ or $\D_n(\tau,\gamma,N)$ as {\sl frequency vectors}, and $\alpha$ in $\D^1_n(\tau,\gamma)$
or $\D^1_n(\tau,\gamma,N)$ as {\sl direction vectors}.

The parameter $N$  is called the {\sl truncation order}, or simply the {\sl cutoff}. We provide more details about
Diophantine sets and the significance of the cutoff below in \S 4.2.  For now, we note that, as discussed in \S 4.2 (i) below, 
$\D^1_n(\tau,\gamma)$ and thus also its supersets are nonempty for $\tau>n-1$ and sufficiently small $\gamma\in(0,1)$.

\section{Basic Problem and Main Result}

Using the notation and terminology above, our basic problem is relatively simple to state.  For fixed $\delta\in(0,\frac{1}{2})$, we seek
the largest set $S\subset\Sno$ of direction vectors $\alpha$ whose corresponding flows $\alpha_t$ fill $\Tn$ quickly, i.e., within a time $T$
depending only (and if possible optimally) on $\delta$ and the parameters determining $S$.

In \S 5.1 below, we give an overview of previous results on this problem.  For now, we repeat that the optimal estimate for the largest set
$S=\D_n(\tau,\gamma,N^*)$ was obtained in 2003 in \cite{BBB03}.
Our contribution in the present paper is to show that geometry-of-numbers methods used by A.~Bounemoura \cite{B16} to get optimal estimates
for $\D_n(\tau,\gamma)$ may be used to recover the results of \cite{BBB03} for the larger set $\D_n(\tau,\gamma,N^*)$.
More specifically, we have the following

\begin{theorem}
For integer dimension $n\ge2$, let $\tau>n-1$, $\gamma\in(0,1)$ be such that $\D^1_n(\tau,\gamma)$ is nonempty.
Choose  $\delta\in(0,\frac{1}{2})$ and set $N^*=(1+n^2n!)/\delta$.
Then given any direction vector $\alpha\in\D^1_n(\tau,\gamma,N^*)$, the flow $\alpha_t$ fills $\Tn$ to within $\delta$
after time $\displaystyle \,T < \frac{C(n,\tau)}{\gamma\delta^\tau}$, \ where $\,C(n,\tau) = (1+n^2n!)^{\tau+1}$.
\end{theorem}

We call $N^*$ the {\sl critical truncation order} or {\sl critical cutoff.}
From the definitions of Diophantine sets, for fixed $n$, $\tau$, $\gamma$, and any $N\ge N^*\ge1$, we have the inclusions
$\D^1_n(\tau,\gamma)\subseteq\D^1_n(\tau,\gamma,N)\subseteq\D^1_n(\tau,\gamma,N^*)$.  
These immediately give the following corollaries of Theorem 1.

\begin{corollary}
For integer dimension $n\ge2$, let $\tau>n-1$, $\gamma\in(0,1)$ be such that $\D^1_n(\tau,\gamma)$ is nonempty.
Then given  $\delta\in(0,\frac{1}{2})$ and any direction vector $\alpha\in\D^1_n(\tau,\gamma)$, the flow $\alpha_t$ fills $\Tn$ to within $\delta$
after time $\displaystyle  \,T < \frac{C(n,\tau)}{\gamma\delta^\tau}$, \ where $\,C(n,\tau) = (1+n^2n!)^{\tau+1}$.
\end{corollary}

\begin{corollary}
For integer dimension $n\ge2$, let $\tau>n-1$, $\gamma\in(0,1)$ be such that $\D^1_n(\tau,\gamma)$ is nonempty.
Then given $\delta\in(0,\frac{1}{2})$, a cutoff $N\ge N^*=(1+n^2n!)/\delta$,
and a direction vector $\alpha\in\D^1_n(\tau,\gamma,N)$, the flow $\alpha_t$ fills $\Tn$ to within $\delta$
after time $\displaystyle  \,T < \frac{C(n,\tau)}{\gamma\delta^\tau}$, \ where $\,C(n,\tau) = (1+n^2n!)^{\tau+1}$.
\end{corollary}

Although these corollaries are simply weaker versions of Theorem 1 (since their hypotheses are more restrictive),
we state them here for their potential use in applications, and for comparison with other results, as discussed further
below in \S 5.2.

\section{Diophantine Sets and the Dynamics of Linear Flow}

In this section, we provide some background for the reader who may be unfamiliar with the connection between
Diophantine sets and the filling property for linear flow on $\Tn$.  The presentation is elementary and 
informal, and many facts are stated without proof.   For further details, we recommend the texts
\cite{BHS96} and \cite{LM88}.

\subsection{Resonance}

To see the connection between Diophantine sets and linear flow, we need some terminology pertaining
to the phenomenon of resonance as it arises in small divisor theory of dynamical systems.  For background and more details,
see Appendix 3 of \cite{LM88}, where what we call resonant multiplicity is instead called resonant dimension, and note that
our situation is simplified substantially by the lack of variables $I$ in an action-like base space (where resonance is usually studied)
and a ``frequency map" $I\mapsto\alpha$ to our frequency vectors $\alpha$.

We say that the frequency vector $\alpha\in\Rn$ is {\sl resonant} if there is a $k\in\Zno$ such
that $k\cdot\alpha=0$.  We denote the set of resonant frequency vectors by $\cR$, and its complement,
the set of {\sl nonresonant} frequency vectors, by $\cN$.  It is not hard to see that both $\cR$ and $\cN$ are
dense in $\Rn$, while $\cR$  is of Lebesgue measure 0 and $\cN$ is of full measure.

Given $k\in\Zno$, the {\sl simple resonance} (or {\sl resonance of multiplicity 1})
determined by $k$ is the hyperplane through the origin $\cR_k=\{\alpha\in\Rn\,\big|\,k\cdot\alpha=0\}$.  
More generally, $\alpha\in\Rn$ belongs to a {\sl resonance of multiplicity} $m\in\{1,\ldots,n-1\}$ if
it belongs to $m$ independent simple resonances $\cR_{k^{(1)}},\ldots,\cR_{k^{(m)}}$, in other words if
$\alpha\in\bigcap_{j=1}^m\cR_{k^{(j)}}$ where $\{k^{(1)}\!,\ldots,k^{(m)}\}$ is linearly independent.
Multiple resonances are nested in the sense that whenever a frequency vector belongs to a resonance of multiplicity $m$,
it also belongs to resonances of lesser multiplicity $l\in\{1,\ldots,m-1\}$.

For any simple resonance $\cR_k$, there are precisely two nonzero integer vectors ($-k'$ and $k'\in{\rm span}\{k\}$)
of smallest norm $\|k'\|$ such that $\cR_k = \cR_{k'}=\cR_{-k'}$; this smallest norm $\|k'\|$ is called the {\sl order}
of the simple resonance.
% The order of the multiple resonance $\bigcap_{j=1}^m\cR_{k^{(j)}}$ is the largest order of the simple resonances $\cR_{k^{(1)}},\ldots,\cR_{k^{(m)}}$.

We connect simple resonances with Diophantine sets as follows.  From the definition of
$D_n(\tau,\gamma)$, we see that for each $k\in\Zno$, the set
$\cZ_k = \{\alpha\in\Rn \big|\ |k\cdot\alpha|<\gamma\|k\|^{-\tau}\}$
is excluded from $D_n(\tau,\gamma)$ (for simplicity, we suppress dependence of $\cZ_k$ on parameters $n$, $\gamma$, $\tau$).
Geometrically, $\cZ_k$ is an open ``hyperslab" centered on $\cR_k$, of half-thickness $\gamma\|k\|^{-\tau-1}$.  In other words,
$\cZ_k$ is the set of points $\alpha$ between the two affine hyperplanes $k\cdot\alpha = \pm  \gamma\|k\|^{-\tau}$.
The thickest such hyperslab containing the simple resonance $\cR_k$ is $\cZ_{k'}$, where $k'$ is one of the two shortest integer vectors
determining $\cR_k$ (i.e., $\|k'\|$ is the order of $\cR_k$).  We call this thickest hyperslab $\cZ_{k'}$ the {\sl resonant zone} around $\cR_k$.

\subsection{The structure of Diophantine and truncated Diophantine sets}

In this subsection, we discuss the parameter values for which the Diophantine sets (truncated or not)
are nonempty, the resonance properties of vectors in them, and the geometry and topology of these sets.

 \medskip
 \noi{\bf (i) \ Diophantine sets are nonempty for $\tau>n-1$ and small $\gamma$}

\noi By writing the set of Diophantine frequency vectors in the form
$\D_n(\tau,\gamma) = \Rn\,\backslash\bigcup_{0\not=k\in\Zn}\cZ_k$,
it becomes a simple exercise to estimate (crudely) its relative Lebesgue measure, as we now outline.
Denoting the closed unit ball in $\Rn$ by $B$ and Lebesgue measure by $\mu$, we readily see that
$
\mu\bigl(B\cap\D_n(\tau,\gamma)\bigr)
\ge
\mu(B)\,- \sum_{0\not=k\in\Zn}\mu(\cZ_k\cap B)  \ge  \mu(B) - \,\gamma \,a_n\sum_{0\not=k\in\Zn}\|k\|^{-\tau-1}
$
where $a_n>0$ is an appropriate constant.
Now the series $\sum_{0\not=k\in\Zn}\|k\|^{-\tau-1}$ converges precisely for $\tau>n-1$; in this case we set
$b(n,\tau) = a_n\sum_{0\not=k\in\Zn}\|k\|^{-\tau-1}$, and we have $\mu\bigl(B\cap\D_n(\tau,\gamma)\bigr)\ge\mu(B) - \gamma \,b(n,\tau)$.
This shows that, for any $\tau>n-1$ and for sufficiently small $\gamma$, the measure $\mu\bigl(B\cap\D_n(\tau,\gamma)\bigr)$ is positive
and  thus $\D_n(\tau,\gamma)$ is nonempty.  We can also see that for $\tau>n-1$, the relative measure of the complement of
$\D_n(\tau,\gamma)$ is $O(\gamma)$ as $\gamma\to0^+$.   A similar argument gives an analogous result for $\D^1_n(\tau,\gamma)$
as a subset of $\Sno$.  Of course, the truncated Diophantine sets are nonempty under the same conditions, since
$\D_n(\tau,\gamma,N)\supset\D_n(\tau,\gamma)$ and $\D^1_n(\tau,\gamma,N)\supset\D^1_n(\tau,\gamma)$.

 \medskip
 \noi{\bf (ii) \ Resonance properties of vectors in Diophantine sets}

\noi Not only do $\D_n(\tau,\gamma)$ and $\D^1_n(\tau,\gamma)$ contain no resonant vectors
(since $\bigcup_{0\not=k\in\Zn}\cR_k \subset \bigcup_{0\not=k\in\Zn}\cZ_k$) but the exclusion of resonant zones $\cZ_{k'}$ around
each resonance $\cR_k$ means that remaining vectors are ``far from resonance,"  or  ``highly nonresonant," 
with the distance of exclusion diminishing with the order of resonance.
By contrast, although $\D_n(\tau,\gamma,N^*)$ and $\D^1_n(\tau,\gamma,N^*)$ maintain these exclusions up to order $N^*$, beyond this
order, no resonances are excluded; the truncated Diophantine sets contain infinitely many vectors resonant at orders higher than $N^*$.

 \medskip
 \noi{\bf (iii) \ Geometry and topology of Diophantine sets}

\noi First, we know that $\D_n(\tau,\gamma)$ is a closed subset of $\Rn$, since its complement $\bigcup_{0\not=k\in\Zn}\cZ_k$
is open.  Second, $\D_n(\tau,\gamma)$ has a simple radial structure: Given any
$\alpha\in\D_n(\tau,\gamma)$ and any $r\ge1$, it follows immediately from the definition of $\D_n(\tau,\gamma)$ that 
$r\alpha\in\D_n(\tau,\gamma)$.  This shows that $\D_n(\tau,\gamma)$ is a collection of
closed half lines directed outward from
the origin in $\Rn$.  The endpoints of the half lines cannot be closer than distance $\gamma$ to the origin, since the thickest
resonant zones ($\cZ_{k'}$ with $\|k'\|=1$) contain the open ball of radius $\gamma$.
Finally, the complement of
$\D_n(\tau,\gamma)$ contains the dense set $\cR$ (the resonant points), so $\D_n(\tau,\gamma)$ has empty interior.  Since
$\D_n(\tau,\gamma)$ is closed with empty interior, it is nowhere dense.

The authors H.K. Broer, G.B. Huitema, and M.B. Sevryuk go further in their description of the Diophantine sets.
In \S 1.5.2 of \cite{BHS96}, they show that the set $\D^1_n(\tau,\gamma) = \D_n(\tau,\gamma)\cap\Sno$
is the union of a countable set and a Cantor set, and they call $\D_n(\tau,\gamma)$ (which they denote by ${\bf R}^n_\gamma$)
a ``Cantor bundle of closed half-lines."  Many authors
refer informally to both $\D^1_n(\tau,\gamma)$ and $\D_n(\tau,\gamma)$ as Cantor sets or Cantor-like sets.

The topology of the truncated Diophantine sets is quite different.
We may write $\D_n(\tau,\gamma,N) = \Rn\,\backslash\bigcup_{0<\|k\|\le N}\cZ_k$, which emphasizes the construction of
$\D_n(\tau,\gamma,N)$ by the removal from $\Rn$ of finitely many hyperslabs $\cZ_k$.  In fact, the relation between $\D_n(\tau,\gamma)$
and $\D_n(\tau,\gamma,N)$ is analogous to the relation between the Cantor ternary set in $\R$ and the finite collection of closed subintervals obtained
at the $N$th step of its construction.  While $\D_n(\tau,\gamma)$ has the complicated topology of a Cantor-like set, its approximating superset 
$\D_n(\tau,\gamma,N)$ has a simple structure: it consists of finitely many closed connected components, each having nonempty interior
and boundary formed by (portions of) hyperplanes.  The sets of direction vectors $\D^1_n(\tau,\gamma)$ and $\D^1_n(\tau,\gamma,N)$
inherit a very similar relationship, since they are simply the intersections of $\D_n(\tau,\gamma)$ and $\D_n(\tau,\gamma,N)$ with $\Sno$.

The difference between the Diophantine and truncated Diophantine sets has significant practical consequences.
In order to decide whether a vector $\alpha$ belongs to the Cantor-like sets
$\D_n(\tau,\gamma)$ or $\D^1_n(\tau,\gamma)$ we must check infinitely many inequalities, in other words we must specify $\alpha$ with
infinite precision.  To see whether $\alpha$ belongs to $\D_n(\tau,\gamma,N)$ or $\D^1_n(\tau,\gamma,N)$, we need check
only finitely many inequalities; this is roughly analogous to determining if a real number belongs to a closed subinterval of $\R$.

A conversation one of us (HSD) had years ago with a theoretical physicist serves to illustrate this last point.  The physicist wished to apply a theorem
from dynamical systems to a mathematical model of particle accelerator dynamics using realistic data, but the hypotheses of the theorem included
Diophantine conditions on the frequency vector.  ``I can't check infinitely many inequalities," said the physicist, to which HSD replied ``You only need
to check them up to a certain order."  The physicist then asked ``But what is that `certain order' precisely?  And can I be sure that the theorem still applies
rigorously when I do that?"  We observe that a result like our Theorem~1 responds positively to the physicist's questions.

\subsection{The dynamics of linear flow on $\Tn$}

Although linear flow on $\Rn$ is very simple, it is not entirely trivial on $\Tn$ because the torus is compact and in some sense ``multiply periodic."
On $\Tn$, there is a basic distinction between the dynamics of linear flow arising from nonresonant versus resonant frequency vectors.
If $\alpha\in\cN$, the linear flow $\alpha_t$ is minimally
ergodic\footnote{``Minimally ergodic" is short for ``minimal and ergodic," where minimal means that every orbit of the flow is dense in $\Tn$.}
on $\Tn$.  If $\alpha\in\cR$, the flow is not ergodic on $\Tn$;
in fact, if $\alpha\in\Rn$ is resonant with multiplicity $m$, then $\alpha_t$ foliates $\Tn$ into invariant ``subtori" of dimension $n-m$,
and $\alpha_t$ is minimally ergodic on each subtorus instead.  

Nevertheless, for fixed $\delta>0$, resonant flow may quickly fill $\Tn$ to within $\delta$, provided $\alpha$ is resonant at high order, because
the subtorus on which $\alpha_t$ is invariant may itself be $\delta$-dense in $\Tn$.
This is most readily seen in the lowest dimension $n=2$, where all resonances are simple and resonant $\alpha$ generate periodic orbits invariant
on subtori of dimension 1 (topological circles). 
To consider specific examples, for $(0,0)\not=(a,b)\in\R^2$ we denote the normalization of $(a,b)$ by $N(a,b)\,$ [i.e., $N(a,b)=(a,b)/\sqrt{a^2+b^2}\,$].
Now let $q\in\Z_+$ and consider the direction vector 
$\alpha=N(q,1)$ which is resonant at order $\sqrt{q^2+1}$.  %\approx q$.
The corresponding flow $\alpha_t$ is not ergodic on $\T^2$, yet fills $\T^2$ to within
$\delta = 1/(2\sqrt{q^2+1})$    %\approx\frac{1}{2q}$
after time $T=\sqrt{q^2+1}$,  %\approx q$, 
which is no doubt the shortest possible filling time for this $\delta$.
On the other hand, ergodicity does not by itself ensure rapid filling.  For any $q\in\Z_+$, the flow $\beta_t$ of the nonresonant direction vector
$\beta = N(q,\sqrt{2})$ is ergodic on $\T^2$, but for fixed $\delta\in(0,\frac{1}{2})$, the filling time becomes arbitrarily long with increasing $q$.

Clearly, the filling property does not coincide with ergodicity on $\Tn$.
Ergodicity is an asymptotic phenomenon realized over infinite time intervals, ensuring that filling occurs for every $\delta>0$.
The filling property, as defined in \S 2 and as used in applications, is realized over finite time intervals for fixed $\delta>0$.
This in a nutshell is why non-ergodic, high-order resonant flow may also fill $\Tn$ quickly, and why
the set $\D^1_n(\tau,\gamma,N^*)$ is a better approximation than $\D^1_n(\tau,\gamma)$ to the largest set $S$ of direction vectors whose flows
quickly fill the torus.  (It is also one reason we now prefer  ``filling time" over the term
``ergodization time" used in our earliest discussions \cite{D88,D91}.)

\section{Previous Results, Optimality, and Applications}

In this section, we place filling-time results in context by briefly discussing their origins and development, their optimality, and their various uses.

\subsection{Previous results on filling times}

Some preliminary remarks are pertinent here.  We point out that the results discussed below
are not always precisely comparable without slight adjustments; this is usually because of variations in the definition of filling (e.g., the filling radius
$\delta$ is sometimes replaced by a filling diameter $\Delta=2\delta$), or in the definition of Diophantine sets (e.g., different norms are used).  We won't detail
these minor variations in the discussion below.  Speaking more broadly, now that filling-time results are relatively mature, we think it's important to highlight
various authors' contributions, especially since a number of advances appear in papers where filling times were used as a tool in the proof of other results,
and so haven't always received the attention they deserve by themselves.

We now briefly describe previous results in chronological order.

\smallskip
\noi{\bf Result (i).} \ Although the filling-time problem is clearly related to earlier results in ergodic theory and uniform distribution, to our knowledge, the first
explicit definitions and treatments as outlined above in \S 2 and \S 3 were given in the thesis \cite{D88} (presenting a mathematical theory of charged particle
motions in crystals) and subsequent article \cite{D91}.
Here Fourier series methods are used first to show that if $\alpha\in\D_n(\tau,\gamma)$, then $\alpha_t$ fills $\Tn$ to within $\delta$ after time
\ $T\le K\gamma^{-1}\delta^{-(\tau+n/2)}$
where $K=K(n,\tau)$ is a suitable constant.  For the truncated Diophantine case, 
it is also shown that there is an $N^\star =  N^\star(\delta,n,\tau)$ such that if
$N> N^\star$ and $\alpha\in\D_n(\tau,\gamma,N)$, then $\alpha_t$ fills $\Tn$ to within $\delta$ after a time which is $T$ (above) multiplied by a messy factor
involving $N$, $N^\star$, $n$, and $\tau$.  These results may be compared (unfavorably) to our Corollaries 1 and 2 in \S 3.
See Theorems 1 and 2 of \cite{D91} for details.

\smallskip
\noi{\bf Result (ii).} \ A few years later, in \cite{CG94}, L.~Chierchia and G.~Gallavotti made use
of a filling-time estimate in their treatment of  ``Arnold diffusion," a kind of instability occurring in Hamiltonian dynamical systems
(see their brief discussion following Eq.~(8.26) on p.~62 of \cite{CG94}).  Though the authors
do not write down the proof in \cite{CG94}, during a later discussion with one of us (HSD),
Gallavotti explained that he had not been aware of previous results, but had used his own simple Fourier series techniques to get the filling time
$T\le T_0\gamma^{-1}\delta^{-(\tau+n)}$ (with suitable $T_0=T_0(\tau,n)$) for $\alpha\in\D_n(\tau,\gamma)$, using only the special case $n=2$ in \cite{CG94}.
We mention this because it seemed remarkable that the problem arose independently in
a different application than Result (i), yet the solution involved Fourier series and gave a similar estimate, namely a power law of the form
$T\sim \delta^{-(\tau+bn)}$ with $b\ge0$.

\smallskip
\noi{\bf Result (iii).} \ Next, in the article \cite{DDG96}, by L.~Dumas, F.~Golse and one of us (HSD), we used filling-time estimates to understand features of
the kinetic theory and mean free path for a Lorentz gas 
in a periodic array of obstacles.  For dimensions $n\ge3$, we used the estimates in \cite{D91}, but for $n=2$, we wrote down
a proof based on continued fractions developed earlier while working on \cite{D88} and mentioned in Remark 3.2 of \cite{D91}.  
We show that if $\alpha\in\D_2(\tau,\gamma)$, then $\alpha_t$ fills $\T^2$ to within $\delta$ after time 
$T \le C' \gamma^{-1}\delta^{-\tau}$ with $C' = 3^\tau\,2^{(\tau+1)/2}$.
The proof also shows that the same filling time holds whenever $N\ge N^*=3/\delta$ and $\alpha\in\D_2(\tau,\gamma,N)$.
As we explain  in \S 5.2 below, these results in the special case $n=2$ are optimal in terms of their power-law dependence on $\delta$.

\smallskip
\noi{\bf Result (iv).} \  In the course of their detailed treatment of the periodic Lorentz gas problem  \cite{BGW98}, J.~Bourgain, F.~Golse, and B.~Wennberg
obtained the filling-time estimate $T \le C'' \gamma^{-1}\delta^{-\tau}$ (with suitable $C''=C''(\tau,n)$) for $\alpha\in\D_n(\tau,\gamma)$, any $n\ge2$,
using Fourier series methods and combinatorial arguments (see Theorem D of  \cite{BGW98}, which may be compared with our Corollary 1 in \S 3).
Again, as explained below, this estimate is optimal in terms of its power-law dependence on $\delta$.

\smallskip
\noi{\bf Result (v).} \  In 2003, M. Berti, L. Biasco, and P. Bolle presented a new approach to the Arnold diffusion problem in which they also developed
their own filling-time estimates implying both Result~{(iv)} and our present Theorem 1.
(Of course the authors' own constant---call it $C'''$---stands in place of our $C$ in Theorem 1 and the $C''$ of Result (iv).)
See Theorems 4.1 and 4.2 of \cite{BBB03} and their proofs in Appendix~B.
We should say that these results are already the best to date in terms of optimality and the size of the set of frequency vectors to which they apply.
We briefly discuss the proof and its relation to our own proof below in \S 6.

\smallskip
\noi{\bf Result (vi).} \   In 2016, A. Bounemoura \cite{B16} recovered the optimal power-law estimate (as above in Results (iv) and (v)) for
$\alpha\in\D_n(\tau,\gamma)$ using geometry-of-numbers methods developed earlier with one of us (SF) in \cite{BF13} and extended further
in the present paper.  In fact, Bounemoura gives the filling time for $\alpha\in\D_n(\tau,\gamma)$ as simply $T\sim\delta^{-\tau}$, but one can
obtain the constants from his Theorem 1, which is more general.  Applying his theorem in the case $\alpha\in\D_n(\tau,\gamma)$ gives 
$T \le C'''' \gamma^{-1}\delta^{-\tau}$ with $C'''' = 2^\tau(n^2n!)^{\tau+1}$.  This is closely related to our Corollary 1 in \S 3, as it should be.

 \subsection{Optimality}
 
In terms of the most important parameter $\delta$, all filling-time estimates so far have the form of a power law
$T \sim \delta^{-(\tau + bn)}$ with $b\ge0$.  As explained in Remark 8 of \cite{DDG96} (where $R$ is used in place of $\delta$),
the use of these estimates in kinetic theory  % of the periodic Lorentz gas
shows that we cannot have $b<0$.  In other words $b\ge0$,
thus the optimal (shortest possible) estimate of this form is $T\sim \delta^{-\tau}$.

To summarize, the optimal $\delta$-dependence was first achieved in the special case $n=2$, both for
$\alpha\in\D^1_2(\tau,\gamma)$ and $\alpha\in\D^1_2(\tau,\gamma,N^*)$ with $N^*=3/\delta$ (Result (iii) above).
It was then obtained for  $\alpha\in\D^1_n(\tau,\gamma)$, any $n\ge2$ (Result (iv)).  Finally, \cite{BBB03} extended the optimal
$\delta$-dependence to $\alpha\in\D^1_n(\tau,\gamma,N^*)$ or $\D_n(\tau,\gamma,N^*)$, any $n\ge2$ (Result (v)).  

We also believe that the $\delta$-dependence of the critical cutoff in the form $N^*\sim\delta^{-1}$ is optimal.
This is because, for certain $\alpha$ resonant at order less than $O(\delta^{-1})$, the flow $\alpha_t$ fills subtori leaving gaps in $\Tn$
larger than $\delta$; in other words the flow fails to fill $\Tn$ to within $\delta$.  A detailed proof of the optimality of $N^*$ along these lines
would involve carefully chosen sequences $\delta_j\to0^+$ and $\{\alpha^{(j)}\}$ 
with $\alpha^{(j)}$ resonant at orders near $\delta_j^{-1}$.

Finally, we say a few words about optimality with respect to parameters other than $\delta$.  
%[ Can we say anything about $\gamma$-dependence, other than ``All estimates have $T\sim\gamma^{-1}$" ??? ]
We do not believe that any of the above constants $C, \,C', C'', C''', C''''$ in \S 5.1 are optimal.  (These depend parametrically on $\tau$ and $n$ and appear
respectively in our Theorem 1 and Results (iii), (iv), (v), (vi) above.)  Without belaboring the point, this is because these constants arise in chains of inequalities
where no special effort was made to ensure sharpness.
If for some reason a sharp or nearly sharp constant of this type is needed, we believe it would be best to estimate it numerically for the desired dimension $n$.
However, we note the following possible future improvement of order constants in our Theorem 1.  As discussed at the end of the appendix below, it is conjectured
that the current bound $n!$ in the Main Duality Result of Successive Minima could be improved to $Kn$ (suitable $K>0$).  If this were achieved, our
constant $C=(1+n^2n!)^{\tau+1}$ could be replaced by $(1+Kn^3)^{\tau+1}$, and our critical cutoff $N^*=(1+n^2n!)/\delta$ by 
$(1+Kn^3)/\delta$.  (Bounemoura's constant $C''''$ could be similarly improved.)

\subsection{Applications}

Since their introduction, filling-time estimates have had a number of applications in dynamical systems and mathematical physics.
Most---but not all---of these applications have connections with nearly integrable Hamiltonian systems (including KAM and/or Nekhoroshev theory),
since Diophantine linear flows on tori occur as a matter of course in such systems.
In this subsection, we list four areas of application in chronological order, citing  a few references and giving a short sketch of each.  We make no
attempt to be comprehensive in our list or references, but rather seek to give the reader a bird's-eye view of how the estimates have been used.

\smallskip
\noi{\bf(i) \ Non-channeling directions in crystals.} \ In their original application, filling-time estimates were used in a mathematical theory describing
the motion of high-energy charged particles as they impinge upon crystals in various directions.  In this setting, direction vectors in a set like 
$\D^1_3(\tau,\gamma,N)$ correspond to so-called non-channeling motions, where a specialized Nekhoroshev theory is used to show that particles undergo
rectilinear motion until they experience close encounters
with crystal nuclei.  The filling-time estimates give an upper bound on the time (and depth) of close encounter, with ramifications in the physics of
ion implantation in crystals.  The use of truncated Diophantine conditions is key here.
 For details, see \cite{D88} or \cite{D93}.

\smallskip
\noi{\bf(ii) \ Arnold diffusion.} \ In the mid 1960s, {V.I.}~Arnold described a mechanism by which a slow, large-scale instability may occur in nearly integrable Hamiltonian systems
with more than two degrees of freedom, even when systems are arbitrarily close to integrable (i.e., even when KAM and Nekhoroshev theorems apply).  This discovery spawned
a large and continuing literature exploring various features of this instability, now loosely called Arnold diffusion.  Arnold's original mechanism uses so-called ``transition chains,"
which include tori supporting Diophantine linear flow.  Unstable orbits stay very near these tori for time intervals determined by filling-time estimates, allowing the (average) speed
of instability to be measured.  This was first done in the previously cited article \cite{CG94} by Chierchia and Gallavotti.  Other researchers, such as {J.-P.}~Marco \cite{M96} and
J.~Cresson \cite{Cr01}, built upon these results with refined techniques and better filling-time estimates.  By contrast, the treatment of Arnold diffusion by Berti, Biasco, and Bolle
\cite{BBB03} avoids the use of transition chains, but uses the best filling-time estimates in a related way to get optimal diffusion times in a particular setting.

\smallskip
\noi{\bf(iii) \ Kinetic theory of the periodic Lorentz gas.} \ This application studies the behavior of a gas of non-interacting point particles moving rectilinearly in an array of
obstacles distributed periodically in space (here ``space" means $\Rn$ with $n\ge2$).  % Interaction of particles with obstacles may include absorption, reflection, or a combination of the two.  
Interest focuses on the behavior of the gas in the so-called macroscopic limit, as the spacing of the array shrinks to zero while the size of obstacles
shrinks at a different rate
controlled by an exponent $\gamma$ (distinct from $\gamma$ used in our Diophantine conditions).  The macroscopic limit depends in turn on the mean free path
(hence on the distribution of free path lengths) of particles in the array.  Filling-time estimates may be used to measure the free path lengths for directions in $\D^1_n(\tau,\gamma)$,
leading ultimately to the existence of a critical exponent $\gamma_c = \frac{n}{n-1}$ dividing gas behavior into three regimes: hydrodynamic behavior for $1\le\gamma<\gamma_c$;
purely ballistic behavior for $\gamma>\gamma_c$; and, most interestingly, behavior possibly governed by a kinetic equation in the so-called Boltzmann-Grad limit $\gamma=\gamma_c$.
This project was begun in \cite{DDG96}, continued in \cite{BGW98}, and has led more recently to a surprising sort of kinetic behavior in the Boltzmann-Grad limit.
We observe that this application of filling-time estimates is not directly connected to nearly integrable Hamiltonian systems.

\smallskip
\noi{\bf(iv) \ Weak KAM theory.} \ This broad and evolving subject lies in the intersection of nonlinear PDE and dynamical systems, using
viscosity solutions of Hamilton-Jacobi equations to find Aubry-Mather sets in Hamiltonian systems.  (Viscosity solutions are a special type of weak solution arising in nonlinear PDE.
Aubry-Mather sets are certain invariant sets of Hamiltonian systems; they include the invariant tori of KAM theory, but also other sets present under weaker conditions.)
Certain Aubry-Mather sets are obtained at points where viscosity solutions have a particular regularity (e.g., Lipschitz continuity or differentiability).  This regularity is
shown, in part, by the use of filling-time estimates.  A procedure of this sort was first carried out by D.A.~Gomes in \cite{Go03}; related techniques have since been used by
K. Wang and J. Yan \cite{WY12},  by K. Soga \cite{So16},  and by H. Mitake and K. Soga \cite{MS18}.  This is an area that might benefit from the use of truncated Diophantine
conditions, as it is connected with specialized numerical methods.

\section{Proofs}

We first state and prove the following proposition, which is a specially adapted version of Proposition~2.3 from \cite{BF13} and is
the most novel part of the present paper. 
We then essentially follow Bounemoura \cite{B16} in using the proposition to prove Theorem 1.

\begin{proposition}
For integer dimension $n\ge2$, choose Diophantine parameters $\tau>n-1$, $\gamma\in(0,1)$
such that $\D^1_n(\tau,\gamma)$ is nonempty, and assume $N>1+n^2n!$.
Let $\alpha\in\D^1_n(\tau,\gamma,N)$.  Then there exist $\omega_1,\ldots,\omega_n\in\Rn$ and $x_1,\ldots,x_n \in\R$ such that, 
for $j\in\{1,\dots,n\}$, we have

\medskip
\noi{\rm (i)}\  \ $\displaystyle \frac{\sqrt{3}}{2} < x_j \le \frac{nn!\, N^\tau}{\gamma}\,$,

\medskip
\noi{\rm (ii)} \ $\displaystyle \|\alpha - \omega_j\| \le \frac{nn!}{x_j(N-1)}\,$, \ and

\medskip\smallskip
\noi{\rm (iii)} \ $\{x_1\omega_1,\ldots,x_n\omega_n\}\,$ is a $\Z$-basis for $\Zn$.
\end{proposition}

In the proof below, we use basic definitions and results from the geometry of numbers.  For the reader's convenience, 
we summarize the needed material in an appendix below.  In what follows, we refer by
[A({\scriptsize R})] to items labeled A({\scriptsize R}) in the appendix, where {\scriptsize R} is a small Roman numeral.

Because Proposition 1 is at the heart of our main result, before beginning the proof, we try to give some insight here into how it works.
The proof of Theorem 1 uses the special $\Z$-basis $\{x_1 \omega_1,\dots, x_n \omega_n\}$ of $\Zn$ from assertion (iii) of Proposition~1.
This basis is separated into multipliers $x_1,\ldots,x_n$ and vectors $\omega_1,\ldots,\omega_n$, which
provides enough flexibility to show two things:
assertion~(i) of Proposition~1, controlling the size of the multipliers as they are used to estimate the filling time, and
assertion~(ii) ensuring that basis vectors $x_j\omega_j$ lie close to the line spanned by $\alpha$, which is used to show that filling occurs.
To achieve this, we first construct a long, slender, solid cylinder $\C$ with central axis $\spalpha$ so that the truncated Diophantine conditions
may be used to check that the reciprocal body $\C^*$ (see [A(i),(iv)]) contains no element of $\Zno$.  
Then $\C,\,\C^*$ may be used in the Main Duality Result of Successive Minima (cf.\ [A(v)]) to show that $nn!\,\C$ contains
the specially adapted $\Z$-basis for $\Zn$ satisfying the required assertions.

\bigskip\noi
{\bf Proof of Proposition 1.} \ Let $H =\bigl(\spalpha\bigr)^\perp$, and consider  $\C,\,\C^*\subset\Rn$ defined by

\noi$\C = \{x\alpha + y\, \big| \,x\in\R, \ y\in H, \ |x|\le N^\tau/\gamma, \ \|y\|\le1/(N-1)\}$ and

\noi$\C^* = \{x\alpha + y \, \big| \, x\in\R, \ y\in H, \ |x|\le \gamma/N^\tau, \ \|y\|\le N-1\}$.

It is a simple matter to verify that $\C,\,\C^*$ are mutually reciprocal CCSBs (compact convex bodies which are symmetric around the origin), as defined below
in the appendix [A(i),(iv)].

We note that if $k\in\Zn\cap\C^*$ then $\|k\| < N$, since $\|k\| = \|x\alpha + y\| \le |x| \|\alpha\|+ \|y\|$, with $|x|\le\gamma/N^\tau<1$, $\|\alpha\|=1$, and $\|y\|\le N-1$.
Now assume $0\not=k\in\Zn\cap\C^*$.  Then combining $\|k\| < N$ with $\alpha\in\D_n(\tau,\gamma,N)$ shows that $|x| = |k\cdot\alpha|\ge\gamma\|k\|^{-\tau}>\gamma /N^\tau$,
which contradicts the definition of $\C^*$.  Therefore $\Zn\cap\C^* = \{0\}$; in other words $\lambda_1(\C^*,\Zn)>1$ by definition
of the first successive minimum [A(ii)].

The main duality result of successive minima [A(v)] reads
$1\le \lambda_k(\C^*,\Zn)\,\lambda_{n+1-k}(\C,\Zn) \le n!$ for $k\in\{1,\ldots,n\}$.
Setting $k=1$ and using $\lambda_1(\C^*,\Zn)>1$ yields $\lambda_n(\C,\Zn) < n!$, which implies [A(i),(iii)] that there is a $\Z$-basis $\{\omb_1,\ldots,\omb_n\}$
of $\Zn$ such that for each $j\in\{1,\ldots,n\}$, \ $\omb_j\in nn!\,\C$.  In other words, $\omb_j = x_j\alpha + y_j$ with
$x_j\in\R$, \ $y_j\in H$, \ $|x_j|\le nn! N^\tau/\gamma$, \ and \ $\|y_j\|\le nn!/(N-1)$.

Since $y_j\in H = \bigl({\rm span}\{\alpha\}\bigr)^\perp\!$, we have $\|\omb_j\|^2 = x_j^2 + \|y_j\|^2$, thus $x_j^2 = \|\omb_j\|^2 - \|y_j\|^2 \ge$ \hfill\break
$1 - \bigl(nn!/(N-1)\bigr)^2 > 3/4$ since we assume $n\ge2$ and $N > 1+ n^2n!$.   Therefore $|x_j|>\sqrt{3}/2$, and changing $\omb_j$ to $-\omb_j$
if necessary, we obtain $x_j >\sqrt{3}/2$.  Together with $|x_j|\le nn! N^\tau/\gamma$, this establishes (i).

Now set $\omega_j = \omb_j/x_j = \alpha + y_j/x_j$, so that $\|\alpha - \omega_j\| = \|y_j/x_j\| \le nn!/x_j(N-1)$, which verifies (ii).
Finally, we see that $\{x_1\omega_1,\ldots,x_n\omega_n\}=\{\omb_1,\ldots,\omb_n\}$ is the $\Z$-basis for $\Zn$ required in (iii).
\hfill$\Box$

\bigskip\noi
{\bf Proof of Theorem 1.} \  Let $\theta\in\Tn$ be arbitrary.  We will prove the theorem by producing a time
$T< (1+n^2n!)^{\tau+1}/(\gamma\delta^\tau)$ such that the endpoint $T\alpha$ of the orbit segment
$\{\alpha_t(0), 0\leq t\leq T\}$ lies within distance $\delta$ of $\theta$.
We use Proposition 1 with $N = N^*= (1+n^2n!)/\delta$.

By Part (iii) of Proposition 1 (and taking into account modular arithmetic on $\Tn$), there exists a unique
$(t_1,\ldots,t_n)\in [0,1)^n$ such that 
$\theta = t_1x_1\omega_1+\cdots+t_nx_n\omega_n$ mod$\,\Zn$.
Set $T =  t_1x_1+\cdots+t_nx_n$.
Then by Part (i) of Proposition 1, we have
$0 \le T =  t_1x_1+\cdots+t_nx_n  \le  x_1+\cdots+x_n  \le  n^2 n! (N^*)^\tau/\gamma = n^2 n!(1+n^2n!)^\tau/(\gamma\delta^\tau)< (1+n^2n!)^{\tau+1}/(\gamma\delta^\tau)$,
as required.
Next, we estimate the distance between $T\alpha$ and $\theta$ as  
$\| T\alpha - \theta\| = \| \sum_{j=1}^n t_j x_j (\alpha - \omega_j)\| \le \sum_{j=1}^n nn!/(N^*-1) = n^2 n!/(N^*-1) < \delta$,
where we use Part (ii) of Proposition 1 in the first inequality, and $N^* = (1+n^2n!)/\delta > 1 + n^2n!/\delta$ in the last inequality.
\hfill$\Box$

\bigskip\medskip\noi
Finally, we say a few words about our proof and its relation to other proofs, especially the proof by Berti, Biasco, and Bolle.
Proofs of the earliest filling-time results made use of Fourier series methods (cf.\ Results (i), (ii), (iv) in \S 5.2 above).  Yet even then,
the short proof of optimal estimates using continued fractions in the special case $n=2$ (Result (iii), \S 5.2) hinted that simpler proofs and better results
would come from number-theoretic methods.  Indeed, though continued fractions don't fully generalize to higher dimensions, a hybrid number-theoretic approach
was found by Berti {\sl et al.}\ in \cite{BBB03}.  We would characterize their approach as a very clever use of geometry-of-numbers methods (without directly using
any of the major theorems of that subject) combined with an induction proof on the dimension. 
By contrast, our approach makes essential use of successive minima in 
geometry of numbers, especially the Main Duality Result of Successive 
Minima (cf.~A(v) in the appendix below).  Using this theorem not only gives a short proof, but we also believe it 
shows filling-time estimates to be a natural part of the geometry of numbers.

\section*{Appendix. \ Geometry of Numbers} 

The geometry of numbers is a branch of number theory begun in the late 19th century by Hermann Minkowski.  It has by now grown substantially
into a vigorous subject in its own right, as the reader may verify by consulting \cite{C59}, \cite{GL87} or other texts.  In this appendix, we provide
only the minimum material needed in the proof of Proposition 1, and we refer to the texts just cited for proofs and further details.

\bigskip\noi
{\bf A(i). Some notation and terminology.}

\noi In the geometry of numbers, it is customary to refer to a connected subset of $\Rn$ with nonempty interior as a ``body."  For simplicity, here
we restrict attention to bodies that are compact, convex, and symmetric around the origin.  We use the abbreviation CCSB to denote such a body.

For $\C\subset\Rn$ and $\lambda\ge0$, we define $\lambda\,\C\subset\Rn$ by $\lambda\,\C = \{\lambda x \,| \,x\in\C\}$, and we sometimes say that
$\lambda\,\C$ is the dilation of $\C$ by $\lambda$, or simply that $\lambda\,\C$ is dilated.

We say that the set $\{\omb_1,\ldots,\omb_n\}\subset\Zn$ is a $\Z$-basis for $\Zn$ if it is linearly independent and if, given any $k\in\Zn$,
there are  $m_1,\ldots,m_n\in\Z$ such that $k=m_1\omb_1+\cdots+m_n\omb_n$.

\bigskip\noi
{\bf A(ii). The $n$ successive minima.}

\noi Given a CCSB $\C\subset\Rn$, for $j\in\{1,\ldots,n\}$ we define the $n$ successive minima
of $\C$ with respect to $\Zn$ by \ 
$\lambda_j(\C,\Zn) = \inf\{\lambda>0 \,|\, {\rm dim\,span}\,(\lambda\,\C\cap\Zn)\ge j\}$.

This says that $\lambda_j$ is the smallest $\lambda$ for which the dilated body $\lambda\,\C$ contains $j$ linearly independent vectors in $\Zn$.
(The definition ordinarily applies when $\Zn$ is replaced by a more general lattice $\Lambda$, but we don't use such $\Lambda$ in this paper.)

\bigskip\noi
{\bf A(iii). Obtaining a $\Z$-basis for $\Zn$.}

\noi Given a CCSB $\C\subset\Rn$, by the definition of $\lambda_n = \lambda_n(\C,\Zn)$, the body
$\lambda_n\,\C$ contains a set of $n$ linearly independent elements of $\Zn$.
This set is not necessarily a $\Z$-basis for $\Zn$, as the $n$-dimensional lattice consisting of its integer combinations may be a proper sublattice of $\Zn$.
However, by dilating further, we can capture a $\Z$-basis:  it follows from the remark after the corollary to Theorem VII in Chapter VIII of \cite{C59}
that the dilated body
$n \lambda_n\,\C$ contains a $\Z$-basis for $\Zn$.

\bigskip\noi
{\bf A(iv). Reciprocal  bodies.}

\noi Given a CCSB $\C\subset\Rn$, we define the corresponding reciprocal body $\C^*\subset\Rn$ by  \hfill\break
$\C^* = \{y\in\Rn\,|\, x\cdot y \le1 \ {\rm for\ all\ } x\in\C\}$.  It is not difficult to show that $\C^*$ is also a CCSB, and that $\C$ is the reciprocal
body corresponding to $\C^*$.  For this reason, we also say that $\C$ and $\C^*$ are mutually reciprocal, or form a mutually reciprocal pair.
(Some authors use the adjectives ``polar" or ``dual" in place of reciprocal, and in \cite{GL87}, the
authors use the compound adjective ``polar reciprocal.")

\bigskip\noi
{\bf A(v). The main duality result of successive minima.}

\noi Given a pair $\C$, $\C^*\subset\Rn$ of mutually  reciprocal CCSBs, for $k\in\{1,\ldots,n\}$ we have   \hfill\break
$1\le \lambda_k(\C^*,\Zn)\,\lambda_{n+1-k}(\C,\Zn) \le n!$. 

This is stated (with $\Zn$ replaced by a general lattice $\Lambda$ and its dual) as Theorem VI of Chapter VIII, \S 5 in \cite{C59}.
A similar theorem, but with upper bound $(n!)^2$ replacing $n!$, is stated as Theorem~5 of Chapter~2, \S14.2 in \cite{GL87}; then in Part ii.6 of the
``Supplement to Chapter 2"  in the same book, it is explained that the upper bound was improved to $n!$ by K.~Mahler in 1939.
In fact, there is an ongoing effort to improve the theorem's upper bound to an optimal value, which is conjectured to be $Kn$ for some
universal constant $K>0$.  See the end of \S 5.2 above for the effect this would have on the order constants in our Theorem 1.

\section*{Acknowledgments}  We thank Abed Bounemoura for helpful discussions, and again for his discovery of the link between
the Main Duality Result (A(v) above) and filling times for linear flow on $\Tn$.

\bibliographystyle{alpha}

\end{document}